\documentclass{article}%
\usepackage{amsfonts}
\usepackage{amssymb}
\usepackage{amsmath}
\usepackage{graphicx}%
\providecommand{\U}[1]{\protect\rule{.1in}{.1in}}

\newenvironment{proof}[1][Proof]{\noindent\textbf{#1.} }{\ \rule{0.5em}{0.5em}}
\begin{document}
\ \ \
\begin{align*}
&  \text{{\LARGE \ Characterisation of the numbers which satisfy}}\\
&  \text{{\LARGE \ \ \ \ \ \ \ \ \ \ \ \ \ \ the\ height reducing property}}%
\end{align*}

\begin{center}
\bigskip

Shigeki AKIYAMA, J\"{o}rg M. THUSWALDNER\ and Toufik ZA\"{I}MI

\bigskip
\end{center}

ABSTRACT.\textbf{\ }\textit{Let }$\alpha$\textit{\ be a complex number. We
show that there is a finite subset }$F$ \textit{of the ring of the rational
integers }$\mathbb{Z},$ \textit{such that} \textit{\ }$F\left[  \alpha\right]
=\mathbb{Z}\left[  \alpha\right]  $\textit{,} \textit{if and only if} $\alpha$
\textit{is an algebraic number whose conjugates, over the field of the
rationals, are all of modulus one, or all of modulus greater than one}.
\textit{This completes the answer to a question, on the numbers satisfying the
height reducing property, posed in Akiyama and Za\"{\i}mi (2013). }

\bigskip

\begin{center}
\textbf{1. Introduction}
\end{center}

Following [1], we say that a complex number $\alpha$ satisfies the height
reducing property, in short HRP, if there is a \ finite subset $F$ of the ring
of the rational integers $\mathbb{Z},$ such that each polynomial with
coefficients in $\mathbb{Z},$ evaluated at $\alpha,$ belongs to the family
$F\left[  \alpha\right]  :=\left\{
{\displaystyle\sum\limits_{j=0}^{n}}
\varepsilon_{j}\alpha^{j}\mid(\varepsilon_{0},...,\varepsilon_{n})\in
F^{n+1},n\in\mathbb{N}\right\}  ,$ where $\mathbb{N}$ is the set of
non-negative rational integers. In this case, we have, by [3, Theorem 1 (i)],
that $\alpha$ is an algebraic number whose conjugates, over the field of the
rationals $\mathbb{Q},$ are all of modulus one, or all of modulus greater than
one (such a number $\alpha$ is called an expanding number [2]). Theorem 1 (ii)
of [3] says also that $\alpha$ satisfies HRP, when it is a root of unity, or
when it is an expanding number. Hence, to obtain a characterisation of numbers
satisfying HRP, it remains to consider the situation where the conjugates of
the algebraic number $\alpha$ belong to the unit circle and are not roots of
unity; this case has been partially treated in [3, Theorem 2], when the
greatest number of multiplicatively independent conjugates of $\alpha,$ over
$\mathbb{Q},$ takes some optimal values.

Recall also, by [2, Theorem 1], that we may suppose that the set $F,$ defined
above, is contained in the complex field $\mathbb{C},$ without affecting the
definition of the HRP; in other words, $\alpha$ satisfies HRP if and only if%

\begin{equation}
\exists\text{\ }F\subset\mathbb{C}\text{ such that \ }\mathbb{Z}\left[
\alpha\right]  =F\left[  \alpha\right]  \text{ \ and }Card(F)<\infty. \tag{1}%
\end{equation}

---------------------------------------------------------------------------------------

\textit{Mathematics Subject Classification} (2010): 11R06, 12D10, 11R04

\textit{Key words and phrases}: Height of polynomials, Special algebraic
numbers, Representations of algebraic numbers

\smallskip

\textit{The second author was supported by projects I1136 and W1230 funded by
the Austrian Science Fund.}

\newpage\ 

By an algebraic approach, we obtain, in this note, that the converse of
Theorem 1 (i) of [3] is true, independently of the distribution, outside the
open unit disc, of the conjugates of $\alpha:$

\medskip\medskip

\textbf{Theorem.} \textit{Let }$\alpha\in\mathbb{C}\mathit{.}$\textit{ Then,
there is a finite subset }$F$ of $\mathbb{Z}$ \textit{such that}
\textit{\ }$F\left[  \alpha\right]  =\mathbb{Z}\left[  \alpha\right]  ,$
\textit{if and only if} $\alpha$ \textit{is an algebraic number whose
conjugates, over }$\mathbb{Q}$\textit{, are all of modulus one, or all of
modulus greater than one}.

\medskip\bigskip

Assume that $\alpha$ satisfies the HRP.  Then, it is natural to ask for the
cardinality of the smallest set $F$ $\subset\mathbb{C}$ (or the smallest set
$F$ $\subset\mathbb{Z})$ satisfying $F[\alpha]=\mathbb{Z}[\alpha].$ This
question is adressed and partially solved in [2]. As mentioned in [3], the
height reducing problem can be compared with canonical number systems and
finiteness property of beta-expansions, where the set $F$ has more specific
shape (some related references may be found in [1, 2, 3, 5]). For example, a
pair $\left(  \alpha,F\right)  ,$ satisfying (1), is called a number system
(resp. a canonical number system) of the ring $\mathbb{Z}\left[
\alpha\right]  ,$ if $0\in F$ and $Card\left(  F\right)  =\left\vert
M_{\alpha}(0)\right\vert $ (resp. if $F=\left\{  0,1,...,\left\vert M_{\alpha
}(0)\right\vert -1\right\}  ),$ where $M_{\alpha}$ designates, throughout, the
minimal polynomial, over $\mathbb{Q},$ of the algebraic number $\alpha$ (the
coefficients of $M_{\alpha}$ are supposed to be rational integers and their
greatest common divisor is one). Recall also that a result of Lagarias and
Wang implies that an expanding integer $\alpha$ satisfies (1) with $F=\left\{
0,\pm1,...,\pm(\left\vert M_{\alpha}(0)\right\vert -1)\right\}  $ [4].

To prove the relation (1) for some fixed pair $\left(  \alpha,F\right)  ,$ it
is generally shown that there exists a positive constant $c=c\left(
\alpha,F\right)  ,$ such that for each $\beta\in\mathbb{Z}\left[
\alpha\right]  ,$ there is some $\varepsilon\in F$ verifying
\[
\frac{\beta-\varepsilon}{\alpha}\in\mathbb{Z}\left[  \alpha\right]  \text{
\ \ and \ \ }\left\Vert \phi_{\alpha}(\frac{\beta-\varepsilon}{\alpha
})\right\Vert <\max\{c,\left\Vert \phi_{\alpha}(\beta)\right\Vert \},
\]
where $\left\Vert .\right\Vert $ is the sup norm (for example) of the
$\mathbb{Q}-$vector space
\[
\mathbb{K}_{\infty}:=\mathbb{R}^{r}\times\mathbb{C}^{s},
\]
$r$ $\left(  \text{resp. }2s\right)  $ denotes the number of real $\left(
\text{resp. of non-real}\right)  $ conjugates, over $\mathbb{Q},$ of the
algebraic number $\alpha,$ and $\Phi_{\infty}$ is the standard Minkowski's
$\mathbb{Q}$-linear map
\[
\Phi_{\infty}:\mathbb{Q}\left(  \alpha\right)  \rightarrow\mathbb{K}_{\infty},
\]
which sends $\alpha$ to its conjugates, over $\mathbb{Q},$ situated in
$\left\{  z\in\mathbb{C}\mid\operatorname{Im}(z)\geq0\right\}  $ (for
example). This allows us sometimes to obtain number systems, when $\alpha$ is
an expanding number, but not when $\left\vert \alpha\right\vert =1$ (see for
instance [2, Section 2]). An alternative solution to this problem is to add
certain finite completions, corresponding to the divisors of the denominator
of the fractional ideal $(\alpha),$ to enlarge the ring $\mathbb{K}_{\infty}$
and the range of the corresponding embedding $\Phi_{\infty}$ : this is the key
of Lemma 1, which is the main result of this manuscript. This lemma is proved
in the last section, and we recall in the next one some related notions.

\smallskip\bigskip

\begin{center}
\textbf{2. Some definitions and notations}
\end{center}

For each given prime $\mathfrak{p}$ of the field~$K:=\mathbb{Q}\left(
\alpha\right)  ,$ where $\alpha$ is a fixed algebraic number, choose an
absolute value $\lvert\cdot\rvert_{\mathfrak{p}}$ in the following way. Let
$\beta\in K$ be given. If $\mathfrak{p}\mid\infty$ \ corresponds to an
Archimedean absolute value, then denote by $\beta^{(\mathfrak{p})}$ the
associated conjugate of~$\beta,$ and set $|\beta|_{\mathfrak{p}}%
=|\beta^{(\mathfrak{p})}|,$ (resp. $|\beta|_{\mathfrak{p}}=|\beta
^{(\mathfrak{p})}|^{2}),$ when $\mathfrak{p}$ is real (resp. is non-real).
With $\mathfrak{p}$ being finite, put $|\beta|_{\mathfrak{p}}=\mathfrak{N}%
(\mathfrak{p})^{-v_{\mathfrak{p}}(\beta)},$ where $\mathfrak{N}(\cdot)$ is the
norm of a (fractional) ideal and $v_{\mathfrak{p}}(\beta)$ denotes the
exponent of $\mathfrak{p}$ in the prime ideal decomposition of the principal
ideal~$\left(  \beta\right)  .$ Write $K_{\mathfrak{p}}$ for the completion of
$K$ w. r. t. \ the absolute value $\lvert\cdot\rvert_{\mathfrak{p}}$ and
recall that this absolute value induces a metric on $K_{\mathfrak{p}}.$

Let $\mathcal{O}$ be the ring of integers of $K,$
\begin{equation}
\alpha\,\mathcal{O}=\frac{\mathfrak{a}}{\mathfrak{b}}\tag{2}%
\end{equation}
where $\mathfrak{a}$ and $\mathfrak{b}$ are coprime ideals in~$\mathcal{O},$
\[
S_{\alpha}=\left\{  \mathfrak{p}:\,\mathfrak{p}\mid\infty
\ \hbox{or}\ \mathfrak{p}\mid\mathfrak{b}\right\}  ,
\]
and define
\[
\mathbb{K}_{\alpha}=\prod_{\mathfrak{p}\in S_{\alpha}}K_{\mathfrak{p}%
}=\mathbb{K}_{\infty}\times\mathbb{K}_{\mathfrak{b}}\,,\quad\mbox{with}\quad
\mathbb{K}_{\infty}=\prod_{\mathfrak{p}\mid\infty}K_{\mathfrak{p}}%
\quad\mbox{and}\quad\ \mathbb{K}_{\mathfrak{b}}=\prod_{\mathfrak{p}%
\mid\mathfrak{b}}K_{\mathfrak{p}}\,.
\]
Then, $\mathbb{K}_{\infty}=\mathbb{R}^{r}\times\mathbb{C}^{s},$ and the
elements of $\mathbb{Q}\left(  \alpha\right)  $ are embedded in $\mathbb{K}%
_{\alpha}$ \textquotedblleft diagonally\textquotedblright\ by the canonical
ring homomorphism
\[
\Phi_{\alpha}:\,\mathbb{Q}\left(  \alpha\right)  \rightarrow\mathbb{K}%
_{\alpha}\,,\quad\beta\mapsto\prod_{\mathfrak{p}\in S_{\alpha}}\beta\,,
\]
where $\ \mathbb{K}_{\alpha}$ is equipped with the product metric of the
metrics induced by the absolute values $\lvert\cdot\rvert_{\mathfrak{p}}.$
Finally, notice that $\mathbb{Q}\left(  \alpha\right)  $ acts multiplicatively
on $\mathbb{K}_{\alpha}$ by the relation
\[
\beta\cdot\left(  z_{\mathfrak{p}}\right)  _{\mathfrak{p}\in S_{\alpha}%
}=\left(  \beta z_{\mathfrak{p}}\right)  _{\mathfrak{p}\in S_{\alpha}},\qquad
\]
where $\beta\in\mathbb{Q}\left(  \alpha\right)  .$

\begin{center}
\textbf{3. Proof of the Theorem}
\end{center}

To make clear the proof of the theorem let us first show three auxiliary
lemmas. The first one is the main tool in this proof.

\medskip

\textbf{Lemma 1.} \textit{Let }$\alpha$ \textit{be} \textit{an algebraic
number, with degree }$n\mathit{,}$\textit{ and } \textit{without conjugates,
over }$\mathbb{Q},$\textit{\ strictly inside the unit circle. Then, there is a
set }$F\subset\mathbb{Z}[\alpha],$ \textit{with cardinality at most}
$2^{n}\left\vert M_{\alpha}(0)\right\vert ,$ \textit{and a constant }%
$c>0$\textit{\ such that for each }$\beta\in\mathbb{Z}\left[  \alpha\right]
,$\textit{\ we can choose }$\varepsilon\in F,$\textit{\ in a way that }%
$\alpha^{-1}(\beta-\varepsilon)\in\mathbb{Z}\left[  \alpha\right]
$\textit{\ with}%

\begin{equation}
\left\vert \alpha^{-1}(\beta-\varepsilon)\right\vert _{\mathfrak{p}}%
<\max\{\left\vert \beta\right\vert _{\mathfrak{p}},c\}, \tag{3}%
\end{equation}
\textit{for each} $\mathfrak{p}\in S_{\alpha}.$

\medskip

\textbf{Proof. }Let $R$ be a complete set of coset representatives of the
finite ring $\mathbb{Z}\left[  \alpha\right]  /\alpha\mathbb{Z}\left[
\alpha\right]  $ and let $\mathcal{U}$ be the collection of the $2^{n}$ open
orthants of $\mathbb{K}_{\infty}\simeq\mathbb{R}^{n},$ where $n=r+2s.$ Since
$\Phi_{\infty}(\alpha\mathbb{Z}\left[  \alpha\right]  )$ contains a lattice
with rank $n$, of $\mathbb{K}_{\infty},$ for each $r\in R$ and each
$U\in\mathcal{U},$ the set $\alpha\mathbb{Z}\left[  \alpha\right]  +r$
contains an element $\varepsilon=\varepsilon\left(  r,U\right)  $ with
$\Phi_{\infty}\left(  \varepsilon\right)  \in U.$ We define the finite set
\[
F=\left\{  \varepsilon\left(  r,U\right)  \;:\;r\in R,\,U\in\mathcal{U}%
\right\}  .
\]
Now, fix $\beta\in\mathbb{Z}[\alpha]$ and pick $\varepsilon\in F$ such that
$\Phi_{\infty}(\varepsilon)$ lies in the same closed orthant as $\Phi_{\infty
}(\beta)$ and satisfies $\alpha^{-1}(\beta-\varepsilon)\in\mathbb{Z}[\alpha].$
It remains to prove that the inequality (3) holds for each $\mathfrak{p}\in
S_{\alpha}.$

Assume first that $\mathfrak{p}\mid\mathfrak{b}.$ Then, as $|\alpha
|_{\mathfrak{p}}>1$ holds by (2), we gain, setting $c_{\varepsilon
,\mathfrak{p}}=|\varepsilon|_{\mathfrak{p}},$ that%

\begin{equation}
|\alpha^{-1}(\beta-\varepsilon)|_{\mathfrak{p}}<\max\left\{  |\beta
|_{\mathfrak{p}},|\varepsilon|_{\mathfrak{p}}\right\}  =\max\left\{
|\beta|_{\mathfrak{p}},c_{\varepsilon,\mathfrak{p}}\right\}  . \tag{4}%
\end{equation}
Next, let $\mathfrak{p}\mid\infty$ be real. Since $\beta^{(\mathfrak{p}%
)}\varepsilon^{(\mathfrak{p})}\geq0$ by the choice of $\varepsilon$ and
$|\alpha|_{\mathfrak{p}}\geq1$ holds by assumption, setting $c_{\varepsilon
,\mathfrak{p}}=2|\varepsilon|_{\mathfrak{p}}$ we have
\begin{equation}
|\alpha^{-1}(\beta-\varepsilon)|_{\mathfrak{p}}\leq|\beta-\varepsilon
|_{\mathfrak{p}}=|\beta^{(\mathfrak{p})}-\varepsilon^{(\mathfrak{p})}%
|<\max\left\{  |\beta|_{\mathfrak{p}},c_{\varepsilon,\mathfrak{p}}\right\}  .
\tag{5}%
\end{equation}
Finally, let $\mathfrak{p}\mid\infty$ be non-real and note that $|\alpha
|_{\mathfrak{p}}\geq1$ holds by assumption also in this case. By the choice of
$\varepsilon$, the complex numbers $\beta^{(\mathfrak{p})}$ and $\varepsilon
^{(\mathfrak{p})}$ lie in the same quadrant of $\mathbb{C}.$ As $\varepsilon
^{(\mathfrak{p})}$ lies in the interior of this quadrant, there is $\eta>0$
depending only on $\varepsilon$ and $\mathfrak{p}$ such that $|\arg{\beta
}^{(\mathfrak{p})}-\arg{\varepsilon}^{(\mathfrak{p})}|<\frac{\pi}{2}-\eta$.
Using this fact, by an easy geometric consideration we obtain%
\begin{equation}
|\alpha^{-1}(\beta-\varepsilon)|_{\mathfrak{p}}\leq|\beta-\varepsilon
|_{\mathfrak{p}}=|\beta^{(\mathfrak{p})}-\varepsilon^{(\mathfrak{p})}%
|^{2}<\max\left\{  |\beta|_{\mathfrak{p}},c_{\varepsilon,\mathfrak{p}%
}\right\}  \tag{6}%
\end{equation}
for some $c_{\varepsilon,\mathfrak{p}}>0$ depending only on $\varepsilon$ and
$\mathfrak{p}$. The inequality (3) now follows from (4), (5) and (6) with
$c=\max\left\{  c_{\varepsilon,\mathfrak{p}}\,:\,\varepsilon\in
F,\,\mathfrak{p}\in S_{\alpha}\right\}  .$%
\endproof

\bigskip\medskip

\textbf{Lemma 2. ([6])} \textit{The ring }$\Phi_{\alpha}\left(  \mathbb{Z}%
[\alpha]\right)  $ \textit{is a discrete subset of} $\mathbb{K}_{\alpha}.$

\medskip

\bigskip\textbf{Proof. }The result is a corollary of Lemmas 3.1 and 3.2 of
[6], where it is shown that $\Phi_{\alpha}\left(  \mathbb{Z}\left[
\alpha\right]  \right)  $ is a Delone set in $\mathbb{K}_{\alpha}.$
\endproof

\bigskip

\textbf{Lemma 3. ([2]) } \textit{If a pair }$\left(  \alpha,F\right)
$\textit{ satisfies the relation }(1)\textit{, then there is a finite subset
}$F^{\prime}$\textit{\ of }$\mathbb{Z}$\textit{ such that }$F^{\prime}\left[
\alpha\right]  =F\left[  \alpha\right]  .$

\medskip

\bigskip\textbf{Proof. }The result follows immediately, by [2, Theorem 1],
where an upper bound (depending only on $\alpha$ and $F)$ of $Card(F^{\prime
})$ is given.%
\endproof

\bigskip

\textbf{Proof of the theorem}. The direct implication is a corollary of
Theorem 1 in [2]. By iterating Lemma~1, we obtain the other implication, using
Lemmas 2 and 3.%
\endproof

\bigskip\bigskip

\bigskip\bigskip\textbf{Acknowledgment.} We thank the referee for careful
reading of this paper.

\medskip

\begin{center}
\bigskip{\LARGE References}
\end{center}

\medskip

[1] S. Akiyama, P. Drungilas and J. Jankauskas, \textit{Height reducing
problem on algebraic integers}, Functiones et Approximatio, Commentarii
Mathematici, \textbf{47} (2012), 105-119.

\medskip

[2] S. Akiyama, J. M. Thuswaldner and T. Za\"{\i}mi, \textit{Comments on the
height reducing problem II, }submitted.

\medskip

[3] S. Akiyama and T. Za\"{\i}mi, \textit{Comments on the height reducing
problem, }Cent. Europ. J. Math., \textbf{11} (2013), 1616-1627.

\medskip

[4] J. C. Lagarias and Y. Wang, \textit{Integral self-affine tiles in}
$\mathbb{R}^{n}.$ \textit{Part II : Lattice tilings}, J. Fourier Anal. Appl.,
\textbf{3} (1997), 83-102.

\medskip

[5] A. Peth\H{o}, \textit{Connections between power integral bases and radix
representations in algebraic number fields}, Proceedings of the 2003 Nagoya
Conference \textquotedblleft Yokoi-Chowla Conjecture and Related
Problems\textquotedblright, 115-125, Saga Univ., Saga, 2004.

\medskip

[6] W.~Steiner and J. M. Thuswaldner, \textit{Rational self-affine tiles,}
Trans. Amer. Math. Soc., to appear.

\medskip

\medskip

\bigskip

Institute of mathematics, University of Tsukuba, 1-1-1 Tennodai, Tsukuba,
Ibaraki, 350-0006 Japan

\textit{E-mail address}: akiyama@math.tsukuba.ac.jp

\bigskip

Department of mathematics and statistics, Leoben University,
Franz-Josef-Strasse 18, A-8700, Leoben, Austria

\textit{E-mail adress:} joerg.thuswaldner@unileoben.ac.at \ \ \ \ \ \ \ \ \ \ \ \ \ \ \ \ \ \ \ \ \ \ \ \ \ \ \ \ \ \ \ \ \ \ \ \ \ \ \ \ \ \ \ \ \ \ \ 

\bigskip

Department of mathematics and informatic, Larbi Ben M'hidi University, Oum El
Bouaghi \ \ 04000, \ Algeria

\textit{E-mail address}: toufikzaimi@yahoo.com\ \ \ 
\end{document}